\documentclass[12pt]{article}
\usepackage{amssymb,latexsym}

\renewcommand{\thefootnote}{}

\begin{document}
\title{Orbits of Certain Endomorphisms
of Nilmanifolds and Hausdorff Dimension}
\author{C.S.Aravinda\\
Chennai Mathematical Institute\\
92 G.N.Chetty Road, Chennai 600 017, INDIA\\
and\\
P.Sankaran\\
Institute of Mathematical Sciences\\
CIT Campus, Chennai 600 113, INDIA \\
E-mail: \begin{tabular}[t]{l}{\tt aravinda@cmi.ac.in}\\
{\tt sankaran@imsc.res.in}  
\end{tabular}\\ [4mm]} 

\date{}

\maketitle

\footnote{2000 A.M.S. Subject Classification:- 53C22, 58F17}
\thispagestyle{empty}

%\input{amssym.def}
%\newsymbol \blacksquare 1004

\def\theequation {\arabic{section}.\arabic{equation}}
\renewcommand{\thefootnote}{}

\newcommand{\codim}{\mbox{{\rm codim}$\,$}}
\newcommand{\stab}{\mbox{{\rm stab}$\,$}}
\newcommand{\lr}{\mbox{$\longrightarrow$}}

\newcommand{\cm}{{\cal M}}
\newcommand{\cf}{{\cal F}}
\newcommand{\cd}{{\cal D}}
%\newsymbol \rightarrowtail 131A
\newcommand{\blr}{\Big \longrightarrow}
\newcommand{\da}{\Big \downarrow}
\newcommand{\ua}{\Big \uparrow}
\newcommand{\hra}{\mbox{\LARGE{$\hookrightarrow$}}}
\newcommand{\rt}{\mbox{\Large{$\rightarrowtail$}}}
\newcommand{\dua}{\begin{array}[t]{c}
\Big\uparrow \\ [-4mm]
\scriptscriptstyle \wedge \end{array}}

\newcommand{\be}{\begin{equation}}
\newcommand{\ee}{\end{equation}}

\newtheorem{guess}{Theorem}%[section]
\newcommand{\bth}{\begin{guess}$\!\!\!${\bf .}~}
\newcommand{\eeth}{\end{guess}}
\renewcommand{\bar}{\overline}
\newtheorem{propo}[guess]{Proposition}
\newcommand{\bpropo}{\begin{propo}$\!\!\!${\bf .}~}
\newcommand{\epropo}{\end{propo}}

\newtheorem{lema}[guess]{Lemma}
\newcommand{\blem}{\begin{lema}$\!\!\!${\bf .}~}
\newcommand{\elem}{\end{lema}}

\newtheorem{defe}[guess]{Definition}
\newcommand{\bdefe}{\begin{defe}$\!\!\!${\bf .}~}
\newcommand{\edefe}{\end{defe}}

\newtheorem{coro}[guess]{Corollary}
\newcommand{\bcor}{\begin{coro}$\!\!\!${\bf .}~}
\newcommand{\ecor}{\end{coro}}

\newtheorem{rema}[guess]{Remark}
\newcommand{\brem}{\begin{rema}$\!\!\!${\bf .}~\rm}
\newcommand{\erem}{\end{rema}}

\newtheorem{exam}[guess]{Example}
\newcommand{\beg}{\begin{exam}$\!\!\!${\bf .}~\rm}
\newcommand{\eeg}{\end{exam}}

\newcommand{\ctext}[1]{\makebox(0,0){#1}}
\setlength{\unitlength}{0.1mm}

\newcommand{\cp}{{\cal P}}
\newcommand{\cu}{{\cal U}}
\newcommand{\bz}{{\Bbb Z}}
\newcommand{\bc}{{\Bbb C}}
\newcommand{\bq}{{\Bbb Q}}
\newcommand{\bt}{{\Bbb T}}
\newcommand{\bh}{{\Bbb H}}
\newcommand{\br}{{\Bbb R}}
\newcommand{\wt}{\widetilde}
\newcommand{\im}{{\rm Im}\,}

\newcommand{\spin}{{\rm Spin}\,}
\newcommand{\ds}{\displaystyle}
\newcommand{\tor}{{\rm Tor}\,}

\def\ns{\mathop{\lr}}
\def\nssup{\mathop{\lr\,sup}}
\def\nsinf{\mathop{\lr\,inf}}
\begin{center}
{\it Dedicated to Professor C.S.Seshadri on the occasion of his 
seventieth birthday}
\end{center}

\noindent
{\bf Abstract:} Let $R$ be an element of $GL(n,\br)$ having integer 
entries and let $\rho:\bt^n \lr \bt^n$ denote the induced map on the 
torus $\bt^n=\br^n/\bz^n$. It is well known that $\rho$ is ergodic with 
respect 
to the Haar measure on $\bt^n$ if and only if none of the eigenvalues of
$R$ is a root of unity. 
Dani has shown that there exists a subset $S$ of $\bt^n$ such that for any 
$x\in S$ and any semisimple surjective endomorphism $\rho$ of $\bt^n$ such that 
the corresponding linear endomorphism has no eigenvalue on the unit circle, the 
closure of the orbit   
$\{\rho^k(x)|k\geq 0\}$ contains no periodic points and that the set $S$ is 
`large' in the sense that 
for any nonempty open set $U$ of $\bt^n$ the set $U\cap S$ has 
Hausdorff dimension $n$. In this paper, we shall prove an analogous result 
for certain endomorphisms of nilmanifolds and infranil manifolds.

\section{Introduction}

Let $N$ be a simply connected nilpotent Lie group. 
Assume that $N$ has a  
discrete subgroup $\Lambda$ so that $M=N/\Lambda$ is compact. (cf. 
Theorem 2.12 and Remark 2.14, \cite{rag}.)
Then $\Lambda\cap N^1$, where $N^1=[N,N]$, is a uniform lattice in $N^1$. 
Let $\bar{\bt}= 
\bar{N}/\bar{\Lambda}\cong N/(N^1\cdot\Lambda)$ where $\bar{N}=N/N^1$ is 
the abelianization of $N$ and $\bar{\Lambda}$ is the 
image of $\Lambda$ under the natural projection $N\lr \bar{N}$.  One has a 
smooth bundle 
projection $p:M \lr \bar{\bt}$ with fibre $N^1/\Lambda\cap N^1$ which 
is again a nilmanifold. 
We put the Riemannian metric on $M$ induced by a right invariant 
metric on $N$.   

Any automorphism $R$ of $N$ such that $R(\Lambda)\subset \Lambda$ 
induces an automorphism $\bar{R}$ of $\bar{N}$ such that 
$\bar{R}(\bar{\Lambda}) \subset \bar{\Lambda}$. $R$ and $\bar{R}$ 
induce self-maps $\rho:M\lr M$ 
and $\bar{\rho}:\bar{\bt}\lr \bar{\bt}$.  The map 
$\rho:M\lr M$ is a bundle map covering the map $\bar{\rho}:\bar{\bt} \lr 
\bar{\bt}$. We call $\rho$ an endomorphism of $M$. 
A theorem of W.Parry \cite{par} says that $\rho$ is ergodic if and only if 
$\bar{\rho}$ is, that is, if and only if none of the 
eigenvalues of $\bar{R}:\bar{N}\lr \bar{N}$ 
is a root of unity.

\bdefe
{\rm Let ${\cal C}$ be a collection of self-maps of a topological space $X$. 
We say that a nonempty subset $S\subset X$ is 
{\it exceptional relative to ${\cal C}$} if the following condition holds:
for any $f\in {\cal C}$ and any $x\in S$, the closure of the $f$-orbit 
$\{f^k(x)|k\geq 0\}$ does not contain any $f$-periodic point}.
\edefe 

When a class ${\cal C}$ of self maps of a smooth Riemannian manifold $M$ 
contains an ergodic map, one expects any exceptional set to be 
`small'. Indeed such a set is necessarily of measure zero and  
is of first category, i.e, its complement 
contains a residual set.  
However, relative to the class of all semisimple surjective  
endomorphisms of the standard $n$-torus $\bt^n$ such that the corresponding 
linear endomorphisms of $\br^n$ have no eigenvalue on the unit circle, Dani 
has proved the existence of an exceptional set $S\subset \bt^n$ which is 
`large' in the sense that for any nonempty open subset $U$ of $\bt^n$, the
set $S\cap U$ has Hausdorff dimension $n$. In fact, Dani shows that $S$ is 
an $\alpha$-winning set of the Schmidt game (cf. \S2).  

The purpose of this note is to extend Dani's result to a certain class of
endomorphisms of nilmanifolds and of infranil manifolds. 
%%%% CHECK 
Examples 
of nilmanifolds are $U/\Gamma$ where $U$ is the group of unipotent 
upper triangular $n$-by-$n$ matrices  over ${\Bbb R}$ (resp. ${\Bbb C}$) 
and $\Gamma$ is the subgroup of those having entries  in ${\Bbb Z}$ (resp. 
in ${\Bbb Z}[i]$).  Nilmanifolds occur naturally in Riemannian 
geometry. For example, they arise in the study of the geometry near  
a cusp in a finite volume 
(non-compact) Riemannian manifold whose sectional curvatures are 
bounded above by a negative constant. See \cite{ebe} for a
detailed study of certain nilmanifolds arising this way.            

Let ${\cal C}_{\bt}$ denote the collection of all surjective 
endomorphisms of an $n$-torus $\bt$ such that the corresponding linear 
automorphisms of $\br^n$ are semisimple and have no eigenvalue on the 
unit circle. We prove 

\bth \label{Haus}
Let $M=N/\Lambda$ be an $n$-dimensional nilmanifold. Let ${\cal C}$ be 
the set of all surjective endomorphisms $\rho$ of $M$ such that 
$\bar{\rho}\in {\cal C}_{\bar{\bt}}$. Then there exists a positive number 
$\mu \leq 1/2$ and an exceptional set $S\subset M$ relative to ${\cal C}$ 
which is $\alpha$-winning for any $\alpha\in (0,\mu]$. In particular,  
$S\cap U$ has Hausdorff dimension $n$ for any nonempty open set $U\subset M$. 
\eeth

Let $N$ be a simply connected nilpotent Lie group. Let $A\!f\!f(N)$ denote 
the group of all affine transformations of $N$ which acts on the right of $N$, 
that is, $A\!f\!f(N)$ is the 
group generated by the right translations by elements of $N$ and the 
automorphisms of $N$. Let $\Gamma$ be a subgroup of $A\!f\!f(N)$ such that 
action of $\Gamma$ on $N$ is free and properly discontinuous with  
the quotient $N/\Gamma$ being compact. The quotient   
$\widehat{M}:=N/\Gamma$ is called an infranil manifold. It is known that  
the subgroup $\Lambda\subset \Gamma$ consisting of all right translations 
of $N$ is a finite index normal subgroup. We shall identify 
an element of $ 
\Lambda$ with the image of the identity element $e\in N$ under it. 
Then it is a result of Auslander \cite{aus} that $\Lambda$ is a 
uniform lattice in $N$.
 
Denote by $M$ the nilmanifold $N/\Lambda$.
One has a finite normal covering $\pi:M\lr \widehat{M}$ with deck 
transformation 
group $G:=\Gamma/\Lambda.$ Starting with the metric on $M$ induced by 
the right invariant 
metric on $N$, the averaging process leads to a $G$-invariant metric 
on $M$. We put this $G$-invariant metric on $M$ and 
the induced metric on $\widehat{M}$ so that $\pi$ is a 
local isometry. 

Let $R:N\lr N$ be an automorphism which preserves the 
$\Gamma$ action, that is, for any $x\in N$ and $\gamma\in \Gamma$, 
there exists a $\gamma'\in \Gamma$ such that $R(x\gamma )=R(x)\gamma'$. 
Then $R(\Lambda)\subset \Lambda$. $R$ induces a self-map 
$\hat{\rho}:\widehat{M}\lr \widehat{M}$ which we call  
an endomorphism of $\widehat{M}$. Note that $\hat{\rho}$ 
is covered by the endomorphism $\rho$ of $M$ defined by $R$. 
We are now ready state

\noindent
\bth \label{infra}   
Let $\widehat{\cal C}$ be the collection of all endomorphisms $\hat{\rho}$
of $\widehat{M}$ induced by automorphisms $R:N\lr N$ such that the 
corresponding toral endomorphism $\bar{\rho}:\bar{\bt}\lr \bar{\bt}$ 
is in ${\cal C}_{\bar{\bt}}$. Then there exist a positive number 
$\nu\leq 1/2$ and an exceptional set $\widehat{S}$ relative  
to $\widehat{\cal C}$ which is $\alpha$-winning for any $\alpha\in (0,\nu]$. 
In particular, for any nonempty open set $U\subset \widehat{M}$, the set 
$\widehat{S} \cap U$ has  Hausdorff dimension $n$.
\eeth

\noindent
{\bf Remark:}
Theorem \ref{infra} includes an important class of manifolds, namely, compact 
connected flat Riemannian manifolds; these are the infranil 
manifolds when the simply connected nilpotent Lie group
$N$ is abelian, i.e., $N = \br^n$. Then the compact flat manifolds are 
precisely those finitely covered by an $n$-torus $\bt^n$.  

Our proofs use Dani's result and some basic observations about 
the nature of pre-image (resp. image) of $\alpha$-winning sets under 
a smooth bundle projection (resp. a covering projection).

\noindent
{\bf Acknowledgements:}
We thank Professor S.G.Dani for his valuable comments 
on an earlier version of this paper. The authors  
gratefully acknowledge financial support from 
Department of Science \& Technology, Government of India.
%We thank the referee of an earlier version of this paper for 
%his/her comments. 

%%%%%%%%%%%%%%%%%%%%%%%  Schmidt  Games  %%%%%%%%%

\section{The Schmidt game}

In this section we give a brief description of the 
Schmidt game \cite{sch} and recall here for the benefit of the 
reader some basic facts about winning sets.

Let $X$ be a complete metric space. Let $0<\alpha, \beta < 1$
and let $S$ be a subset of $X$.
The Schmidt game on $X$ is played 
by two players ${\cal A}$ and ${\cal B}$ as follows: ${\cal B}$ first 
chooses any closed ball $B_0$ in $X$ of radius $r$ where $r\in \br$ is 
positive. Then 
${\cal A}$ 
chooses a closed ball $A_1\subset B_0$ of radius $\alpha r$. 
Then ${\cal B}$ chooses 
a closed ball $B_1\subset A_1$ of radius $r_1= \beta\alpha r$.  
Inductively, after ${\cal B}$ has chosen closed a ball $B_{k}$ of 
radius $r_k=(\alpha\beta)^kr$, $k\geq 1$, ${\cal A}$ chooses a 
closed ball $A_{k+1}\subset B_k$ 
of radius $\alpha r_k$ and ${\cal B}$ chooses a closed ball $B_{k+1}$ of 
radius $r_{k+1}=\beta\alpha r_k$. Since $\displaystyle{\lim}_{k\rightarrow 
\infty}r_k =0$, and since $X$ is a complete metric space, $\cap_{k\geq 1} A_k$ 
is a singleton set $\{x\}$. We say that $S$ is an $(\alpha, \beta)$-winning 
set (for ${\cal A}$) 
if for any choices of ${\cal B}$, ${\cal A}$ can always make her 
choices so that $x\in S$. We say that $S$ is $\alpha$-winning if 
it is $(\alpha,\beta)$-winning for any $\beta\in (0,1)$. 
Of course, the whole space $X$ is 
an $\alpha$-winning set for any $\alpha \in (0,1)$. Schmidt 
has shown that if $\beta\leq 2-\alpha^{-1}$, then the only 
$(\alpha,\beta)$-winning set is $X$ itself.  However, if 
$\beta>2-\alpha^{-1}$, then there can be proper subsets of 
$X$ which are $(\alpha,\beta)$-winning. Schmidt introduced these notions 
in his study of badly approximable numbers and related Diophantine
problems.  Schmidt games have proved to be very useful in 
showing that certain sets which naturally arise in dynamical  
systems and which are small in the 
conventional sense  are indeed large in the sense of Hausdorff 
dimension.  

We state without proofs some basic  
properties of $\alpha$-winning sets: 

\bpropo\label{basic}
\noindent
(i) Suppose $S,~S'$ are $\alpha$-winning sets in $X, X'$ then  
$S\times S'\subset X\times X'$ is also an $\alpha$-winning set.\\ 
(ii) Let  $X_i, ~1\leq i\leq k,$ be closed subspaces of a compact 
metric space $X$ whose interiors cover $X$. If $S_i\subset X_i$ 
is $\alpha_i$-winning in $X_i$ then $S=\bigcup_{1\leq i\leq k}S_i$ 
is $\alpha$ winning where $\alpha =\min\{\alpha_i\mid 1\leq i\leq k\}$.  \\
(iii) {\rm (Prop. 3.3\cite{ara}, cf. \S11,\cite{sch})} Any $\alpha$-winning 
set $S$ in a complete Riemannian manifold $M$ is large, i.e.,  
for any nonempty open set $U\subset M$ the Hausdorff dimension of $S\cap U$ 
is $n$ where $n=\dim M$.\\
(iv) {\rm (\cite{dani1}, cf. Th. 1, \cite{sch})} Let  $f:X\lr Y$ be a 
homeomorphism between 
complete metric spaces and $\lambda$ be a bi-Lipschitz constant  for $f$, that 
is, $\lambda\geq 1$ is a real number such that 
for any $x,y\in X$, one has $\lambda^{-1} d_X(x,y)\leq d_Y(f(x),f(y))\leq 
\lambda d_X(x,y)$. If $S\subset X$ is $(\alpha,\beta)$-winning set in $X$ 
with $\beta\lambda^2<1$,  
then $f(S)$ is an $(\alpha/ \lambda^2,\beta\lambda^{2})$-winning set in $Y$.
In particular, if $S$ is $\alpha$-winning, then $f(S)$ is 
$\alpha/\lambda^2$-winning. \hfill $\Box$
\epropo

\brem \label{nometric}
If $X,Y$ are smooth compact Riemannian manifolds 
and $f:X\lr Y$ is a diffeomorphism, then a bi-Lipschitz constant for 
$f$ always exists (see Ch. 5 of \cite{dani1}, Prop. 3.3 of 
\cite{ara}). Thus, in view of Prop. \ref{basic}$(iv)$,  
the `largeness' of a subset of a compact 
Riemannian manifold does not depend on the specific choice of the metric.  
\erem
%As an immediate corollary of Dani's theorem we obtain the following 

\brem \label{tor}
Let $\Lambda$ be any lattice in the Euclidean space $\br^n$ and 
let $\bt=\br^n/\Lambda$. Consider the collection ${\cal C}_{\bt}$.
Then there exists a subset $S\subset \bt$ which is exceptional 
relative to ${\cal C}_{\bt}$ such that $S$ is $\alpha$-winning for any 
$\alpha\in(0,1/2]$. In particular, $S\cap U$ has Hausdorff dimension 
$n$ for any nonempty open set $U\subset \bt$. 
To see this, note that if $\Lambda=\bz^n\subset\br^n$, then this is just a 
restatement of Dani's result \cite{dan}. In the general case, we observe that 
Dani's proof still goes through when $\bz^n$ is replaced by an arbitrary 
uniform lattice $\Lambda$. 
We can also see this easily using Dani's result for the torus $\br^n/\bz^n$ 
and Remark \ref{nometric}. 
\erem

We conclude this section with the following 

\blem\label{fib}
Let $p:M\lr M'$ be a smooth fibre bundle  
where $M$ and $M'$ are  compact Riemannian manifolds.  
Assume that $S'\subset M'$ is $\alpha$-winning 
for any $\alpha\in (0,\nu]$ for some $\nu\in (0,1/2]$. Then there exists a 
$\mu>0$ such that 
$S:=p^{-1}(S')$ is $\alpha$-winning for any $\alpha\in (0,\mu]$.
\elem

\noindent
{\bf Proof:} Let $F$ be the fibre of the bundle, endowed with 
a Riemannian metric.  
Let $V\subset M'$ be any trivializing open set which is 
diffeomorphic to the  unit ball $B$ in $\br^n$.  
Let $K$ be the image in $V$ of the closed disk of radius $1/2$ 
under the diffeomorphism $B\cong V$. 
Let $h_V:p^{-1}(V)\lr V\times F$ be a diffeomorphism which is 
a trivialization. Then $S'_K:=S'\cap K $ is an $\alpha$- winning 
set in $K$ for any $\alpha\in (0,\nu]$ and so, by Prop. \ref{basic}$(i)$, 
$S'_K\times F$ is an 
$\alpha$-winning set in $K\times F$ for any $\alpha\in (0,\nu]$. Now, 
it follows from Prop. \ref{basic}$(iv)$ that 
$p^{-1}(S'_K)=h_V^{-1}(S_K'\times F)$ is $\alpha$-winning in 
$p^{-1}(K)$
for any $\alpha\in (0,\nu/\lambda^2]$ where $\lambda$ is a bi-Lipschitz 
constant for the diffeomorphism $h_V$.  One can cover $M'$ by the interiors 
of (closed) disk neighbourhoods $K_i,1\leq i\leq k$, where  
$K_i\subset V_i\subset M'$, 
the $V_i$ being trivializing open sets for the bundle $p:M\lr M'$. 
Choose bi-Lipschitz constants $\lambda_i$ for trivializations $p^{-1}(V_i)\lr 
V_i\times F$.  Then $p^{-1}(S'_{K_i})$ is $\alpha$-winning in $p^{-1}(K_i)$
for any $\alpha\in (0,\nu/\lambda_i^2]$ and for each $i$ such that 
$1 \leq i \leq k$. Hence, by Prop. \ref{basic}$(ii)$, it follows that 
$S=\bigcup_{1\leq i\leq k}
p^{-1}(S'_{K_i})$ is $\alpha$-winning for any $\alpha\in (0,\mu]$
where $\mu$ is the smallest of the numbers $\nu/\lambda_i^2, ~1\leq i\leq k$. 
\hfill $\Box$

%%%%%%%%%%%%%%%%%%   PROOF %%%%%%%%%%%%%%

\section{Proofs of Main Theorems}

Let $f:X\lr X$ be a self-map and let $x\in X$. Denote by 
$C_f(x)$ the closure of the set $\{f^k(x)\mid k\geq 0\}$. 

\blem\label{pull}
Let $p:X\lr X'$ be any continuous map. Let ${\cal C}$ and ${\cal C'}$ 
be any collections of self-maps of $X$ and $X'$ respectively such that
given $f\in {\cal C}$, there exists an $f'\in {\cal C'}$ such that 
$p\circ f=f'\circ p$. Suppose that 
$S'\subset X'$ is exceptional relative to ${\cal C'}$, then $S:=p^{-1}(S')$ 
is  exceptional relative to ${\cal C}$.
\elem

\noindent
{\bf Proof:} Let $x\in S,~f\in {\cal C}$. Suppose that $y\in C_f(x)$ is 
$f$-periodic, say, $f^k(y)=y$. Let $x'=p(x)\in S'$ and let $y'=p(y)$. 
Choose $f'\in {\cal C'}$ such that $p\circ f=f'\circ  p$.  Then, 
$p(C_f(x))\subset C_{f'}(x')$. In particular $y'\in C_{f'}(x')$.
Note that $f'^k(y')=f'^k\circ p(y)=p\circ f^k(y)=p(y)=y'$ which
contradicts our hypothesis that $S'$ is exceptional. 
This shows that $S$ is exceptional relative to ${\cal C}$. \hfill $\Box$

\blem \label{push} 
Let $p:X\lr \widehat{X}$ be a finite covering projection where 
$X,~\widehat{X}$  are compact metric spaces. 
Let ${\cal C}$ and $\widehat{{\cal C}}$ be any collections of 
self-maps of $X$ and $\widehat{X}$ respectively such that 
for any $\hat{f}\in \widehat{{\cal C}}$, there exists an $f\in {\cal C}$ 
such that $p\circ f=\hat{f}\circ p$. Let $S$ be an exceptional 
set relative to ${\cal C}$. Then $\widehat{S}:=p(S)\subset \widehat{X}$ is an 
exceptional set relative to $\widehat{{\cal C}}$.
\elem

\noindent
{\bf Proof:} Let $\hat{x}\in \widehat{S}, ~\hat{f}\in \widehat{{\cal C}}$ 
and let $\hat{y}\in C_{\hat{f}}({\hat{x}})$. 
We shall assume that $\hat{y}$ is an $\hat{f}$-periodic point and arrive 
at a contradiction.

Since $\hat{y}\in C_{\hat{f}}(\hat{x})$, there exists 
a sequence of positive integers $n_k$ such that $\hat{f}^{n_k}(\hat{x})$ 
converges to $\hat{y}$. 
Choose $x\in S, ~f\in {\cal C}$ such that $p(x)=\hat{x}$ and 
$p\circ f=\hat{f}\circ p$. 
Consider the set $\{f^{n_k}(x)| k\geq 1\}\subset 
C_f(x)$. Because $x\in S$, this set is infinite. 
Since $X$ is compact, it has a subsequence 
$f^{n_{k_i}}(x)$ converging to an element $y\in X$ and $y\in C_f(x)$ 
since $C_f(x)$ is closed.
The element $y$ must be in the fibre over $\hat{y}\in \widehat{X}$ by 
continuity of $p$. Since $\hat{y}$ is $\hat{f}$-periodic, and since 
$f$ covers $\hat{f}$, it follows that $y$ must be 
$f$-periodic. That is, $y\in C_f(x)$ is an $f$-periodic point which is 
a contradiction since $x\in S,~f\in{\cal C}$ and $S$ is an 
exceptional set relative to ${\cal C}$. \hfill $\Box$

We are now ready to prove the main theorems stated in the introduction. 
We keep the notations of \S1.

\noindent
{\bf Proof of Theorem \ref{Haus}:}  Consider the fibration 
$p: M\lr M'$ where $M'$ is the torus $\bar{N}/\bar{\Lambda}$ and let   
${\cal C}'={\cal C}_{\bar{\bt}}$. By Remark \ref{tor}, there exists 
%number $\nu\leq 1/2$ 
an exceptional set $S'\subset M'$ relative to ${\cal C'}$ which is 
$\alpha$-winning for any $\alpha\in (0,1/2]$. 
Note that for the collection ${\cal C}$ of the theorem, 
the hypotheses of Lemma \ref{pull} are satisfied. It follows that 
$S=p^{-1}(S')$ is exceptional.  By Lemma \ref{fib}, we see that 
$S$ is $\alpha$-winning for any $\alpha\in (0,\mu]$ for some $\mu\in 
(0,1/2]$. \hfill $\Box$ 

\noindent
{\bf Proof of Theorem \ref{infra}:} Consider the finite covering projection 
$\pi:M\lr \widehat{M}$. Recall that the 
deck transformation group $G=\Gamma/\Lambda$ 
acts on $M$ by isometries. 

The hypotheses of Lemma 
\ref{push} holds when we take ${\cal C}$ to be as in Theorem 
\ref{Haus}, that is, the set of all surjective 
endomorphisms $\rho$ of $M$ such that 
$\bar{\rho}\in {\cal C}_{\bar{\bt}}$. 

By Theorem \ref{Haus} there exists a positive 
number $\mu\leq 1/2$ 
and a set $S$ which is  
exceptional relative to ${\cal C}$ such that $S$   
is $\alpha$-winning for any $\alpha\in (0,\mu]$ when $M$ is given the 
right invariant metric. By Prop. \ref{basic}$(iv)$ applied to 
the identity map of $M$, the same statement 
holds if we replace $\mu$ by a suitably smaller positive number $\nu$ when 
we put the $G$-invariant metric on $M$.  
By Lemma \ref{push}, $\widehat{S}:=\pi(S)$ is an exceptional set relative 
to $\widehat{{\cal C}}$. 

Since $G$ acts on $M$ by isometries, the covering projection $\pi$ is a 
local isometry. Choose a  
finite collection $\{K\}$ of normal (closed) disk neighbourhoods 
whose interiors cover 
$M$ and such that $\pi|K:K\lr \pi(K)$ is 
an isometry for each $K$. Then, $S\cap K$ is $\alpha$-winning in $K$ and 
$\pi(S\cap K)$ is $\alpha$-winning for each $\alpha\in (0,\nu]$
by Theorem 1, \cite{sch}.  
It follows from Prop. \ref{basic}$(ii)$ that    
$\widehat{S}=\pi(S)$ is an $\alpha$-winning set for any 
$\alpha\in (0,\nu]$. This completes the proof. \hfill $\Box$

\noindent
{\bf Concluding Remarks:}  
$(i)$  It would be interesting to know whether there are 
large exceptional sets for the class consisting of an 
arbitrary  expanding endomorphism  
(cf. \cite{shu}) of a nilmanifold or an infranil manifold. 

$(ii)$   Kleinbock \cite{kl} has obtained a generalization of Dani's 
theorem showing existence of large sets in a very general 
setting. 
However, Kleinbock's result 
still does not subsume our results because  
the exceptional set in his set up is not known to be an 
$\alpha$-winning set.  This is a stronger conclusion in our 
case which enables us to deduce the maximality of Hausdorff dimension 
for exceptional sets of a {\it family}  of endomorphisms and not just 
for exceptional sets of a single endomorphism.

%%%%%%%%%%%%%%%%%%   References   %%%%%%%%%%%%%%%%

\end{document}